\documentclass[global]{svjour}
\usepackage{amssymb}
\usepackage{amsfonts}
\usepackage{amsmath}

\input{tcilatex}

\begin{document}

\title{{\normalsize Existence of the equilibrium in choice}}
\author{Monica Patriche}
\institute{University of Bucharest 
\email{monica.patriche@yahoo.com}%
}
\mail{University of Bucharest, Faculty of Mathematics and Computer Science, 14
Academiei Street, 010014 Bucharest, Romania}
\maketitle

\begin{abstract}
In this paper, we prove the existence of the equilibrium in choice for games
in choice form. These games have recently been introduced by A. Stefanescu,
M. Ferrara and M. V. Stefanescu. Our results link the recent research to the
older approaches, regarding games in normal form or qualitative games.
\end{abstract}


\section{Introduction}

It has been widely agreed upon the fact that Nash's concept of equilibrium
reflects the possibility of challenging the choices of the unilateral acts
of the players involved in noncooperative games. Having the aim of gaining
the most they can do, the players are in the situation of making choices in
a process described, mathematically, by the notion of "game" (this was
proposed by Nash, and it was initially called "game in normal form"). The
original definitions of Nash \cite{80},\cite{79} have been extended, but the
derived notions of qualitative game or abstract economy and their
corresponding concepts of equilibrium reflect the initial meaning of
flexible elections in any given context to permit the players not to
deviate, once they agreed on the best solution for them.

The new extension of a game, called "the game in choice form", is due to
Stefanescu, Ferrara and Stefanescu \cite{SS}. The game in choice form is a
family of the sets of individual strategies and choice profiles. The authors
also defined the concept of equilibrium in choice for this type of game.
Their interpretations evolved along with the problem of noncooperative
solutions of the voting operators. Firstly, Stefanescu and Ferrara proposed
the concept of Nash equilibrium in choice in \cite{st} and they renamed it
in \cite{SS}, founding conditions under which the equilibrium in choice
exists. The new adopted definitions are coherent with the old underlying
formalism. For instance, when the utility functions represent the players'
options, a choice profile can be seen as the family of players' graphs best
reply mappings. In this case, the set of equilibria in choice coincides with
the set of Nash equilibria. The generality of the new concepts raise the
interest of the formalist theorist of games, which explicitly can show their
significance. The authors themselves developed these ideas in their work.
They referred to the fact that the players' preferences need not be
explicitly represented, at the same time considering the possibility of
recuperating the known solutions as particular cases. The second problem
raised is the one of the nonexistence of a best reply. This interest is
obviously at the core of our original research. In this paper, we are
looking for new conditions, in order to obtain the existence of the
equilibrium in choice. Our assumptions are different from the ones proposed
when the new theory was framed. They concern the properties of the sets of
choice profiles. We are now exploring a method of proof based on continuous
selection and fixed point theorems for correspondences defined by using the
upper sections of the sets which form the game.

The rest of the paper is organised as follows: In Section 2, some notation,
terminological convention, basic definition and results about
correspondences and games in choice form are given. Section 3 contains
existence results for equlibrium in choice.

\section{Preliminaries and notation}

Throughout this paper, we shall use the following notations and definitions:

Let $A$ be a subset of a vector space $X$, $2^{A}$ denotes the family of all
subsets of $A$ and co$A$ denotes the convex hull of $A$. If $A$ is a subset
of a topological space $X,$ cl$A$ denotes the closure of $A$ in $X$. If $T$, 
$G:X\rightarrow 2^{Y}$ are correspondences, then co$G$ and $G\cap
T:X\rightarrow 2^{Y}$ are correspondences defined by $($co$G)(x):=$co$G(x)$
and $(G\cap T)(x):=G(x)\cap T(x),$ for each $x\in X$, respectively.$\medskip 
$

Given a correspondence $T:X\rightarrow 2^{Y},$ for each $x\in X,$ the set $%
T(x)$ is called the \textit{upper section} of $T$ at $x.$ For each $y\in Y,$
the set $T^{-1}(y):=\{x\in X:y\in T(x)\}$ is called the \textit{lower section%
} of $T$ at $y$. The correspondence $T^{-1}:Y\rightarrow 2^{X},$ defined by $%
T^{-1}(y)=\{x\in X:y\in T(x)\}$ for $y\in Y$, is called the \textit{(lower)
inverse} of $T.$

Let $X$, $Y$ be topological spaces. The correspondence $T:X\rightarrow 2^{Y}$
is called \textit{lower semicontinuous} if for each x$\in X$ and for each
open set $V$ in $Y$ with $T(x)\cap V\neq \emptyset $, there exists an open
neighborhood $U$ of $x$ in $X$ so that $T(y)\cap V\neq \emptyset $ for each $%
y\in U$.$\medskip $

The following lemma will be useful to our study of existence of equilibrium
in choice.

\begin{lemma}
(see Yuan \cite{134}). \textit{Let }$X$\textit{\ and }$Y$\textit{\ be two
topological spaces and let }$W$\textit{\ be an open (resp. closed) subset of 
}$X.$\textit{\ Suppose }$T_{1}:X\rightarrow 2^{Y}$\textit{\ , }$%
T_{2}:X\rightarrow 2^{Y}$\textit{\ are upper semicontinuous (resp. lower
semicontinuous) correspondences such that }$T_{2}(x)\subset T_{1}(x)$\textit{%
\ for all }$x\in W.$\textit{\ Then the correspondence }$T:X\rightarrow 2^{Y}$%
\textit{\ defined by}
\end{lemma}

\begin{center}
$T\mathit{(x)=}\left\{ 
\begin{array}{c}
T_{1}(x)\text{, \ \ \ \ \ \ \ if }x\notin W\text{, } \\ 
T_{2}(x)\text{, \ \ \ \ \ \ if }x\in W%
\end{array}%
\right. $
\end{center}

\textit{is also upper semicontinuous (resp. lower semicontinuous).\medskip }

Further, we present the main models of noncooperative games we will deal
with in this paper. The corresponding notions of equilibrium are also
recalled.

Let $(X_{i})_{i\in N}$ be the family of the individual sets of strategies
and let $X=\prod\nolimits_{i\in I}X_{i}$.

The \textit{normal form} of an $n$-person game is $(X_{i},r_{i})_{i\in N}$,
where, for each $i\in N$, $X_{i}$ is a nonempty set (the set of individual
strategies of player $i$) and $r_{i}$ is the preference relation on $X$ of
player $i$. The individual preferences $r_{i}$ are often represented by 
\textit{utility functions}, i.e. for each $i\in N$ there exists a real
valued function $u_{i}:X\rightarrow \mathbb{R}$ (called the utility function
of $i$), such that: $xr_{i}y\Leftrightarrow u_{i}(x)\geq u_{i}(y),$ $\forall
x,y\in X.$ Then the \textit{normal form of n-person game} is $%
(X_{i},u_{i})_{i\in N}$.

We\textit{\ }denote $x_{-i}=(x_{1},...,x_{i-1},x_{i+1},...,x_{n})$, $%
X_{-i}=\prod\nolimits_{i\neq j}X_{i}$ and $(x_{-i},X_{i})=%
\{(x_{-i},x_{i}):x_{i}\in X_{i}\}.$

The\textit{\ Nash} \textit{equilibrium }for the game $(X_{i},u_{i})_{i\in N}$
is a point $x^{\ast }\in $ $X$ which satisfies for each $i\in N:$ $%
u_{i}(x^{\ast })\geq u_{i}(x_{-i}^{\ast },x_{i})$ for each $x_{i}\in X_{i}.$

For each $i\in N,$ the player's i's best reply mapping is the correspondence 
$B_{i}:X_{-i}\rightarrow 2^{X_{i}},$ defined by $B_{i}(x_{-i})=\{x_{i}\in
X_{i}:u_{i}(x_{-i},x_{i})\geq u_{i}(x_{-i},y_{i})$ for each $y_{i}\in
X_{i}\}.$ Then, $x^{\ast }\in $ $X$ is a Nash equilibrium if\textit{\ }only
if\textit{\ }$x^{\ast }\in \bigcap\nolimits_{i\in N}$Gr$(B_{i}).$

The element $x^{\ast }\in $ $X$ is called \textit{weak Nash equilibrium}
(Stefanescu, Ferrara, Stefanescu \cite{SS}):

$u_{i}(x^{\ast })\geq u_{i}(x_{-i}^{\ast },x_{i})$ for each $x_{i}\in X_{i}$
and $i\in N$ such that $B_{i}(x_{-i}^{\ast })\neq \emptyset .\medskip $

A \textit{qualitative game\ }$\Gamma =(X_{i},P_{i})_{i\in N}$ is defined as
a family of $n$ ordered triplets $(X_{i},P_{i}),$ where for each $i\in N$: $%
P_{i}:X\rightarrow 2^{X_{i}}$ is a preference\ correspondence. An \textit{%
equilibrium} for $\Gamma $ is a point $x^{\ast }\in $ $X$ which satisfies
for each $i\in N:$ $P_{i}(x^{\ast })=\emptyset .$ A \textit{weak equilibrium}
(Stefanescu, Ferrara, Stefanescu \cite{SS}) of $\Gamma =(X_{i},P_{i})_{i\in
N}$ is a point $x^{\ast }\in $ $X$ which satisfies $P_{i}(x^{\ast
})=\emptyset $ for each $i\in N$ such that $\{x_{i}\in
X_{i}:P_{i}(x_{-i}^{\ast },x_{i})=\emptyset \}\neq \emptyset .$

A \textit{choice profile} (Stefanescu, Ferrara, Stefanescu \cite{SS}) is any
collection $C:=(C_{i})_{i\in N}$ of nonempty subsets of $X.$ A \textit{game
in choice form} is a double family $((X_{i})_{i\in N},(C_{i})_{i\in N})$,
where $C=(C_{i})_{i\in N\text{ }}$is a choice profile.

We denote $C(x_{-i})$ the upper section through $x_{-i}$ of the set $C_{i},$
i.e., $C(x_{-i})=\{y_{i}\in X_{i}:(x_{-i},y_{i})\in C_{i}\}.$

We will make the following assumption:

(A) assume that for each $x\in X,$ there exists $i\in N$ such that $%
C_{i}(x_{-i})\neq \emptyset .$

The game strategy $x^{\ast }\in X$ is an \textit{equilibrium in choice}
(denoted EC) (A. Stefanescu, M. Ferrara, V. Stefanescu \cite{SS}) if $%
\forall i\in N,$ $(x_{-i}^{\ast },X_{i})\cap C_{i}\neq \emptyset \Rightarrow
x^{\ast }\in C_{i},$ equivalently, $x_{i}^{\ast }\in C_{i}(x_{-i}^{\ast }),$
for every $i\in N$ for which $C_{i}(x_{-i}^{\ast })\neq \emptyset .$ The
strategy $x^{\ast }\in X$ is a \textit{strong} \textit{equilibrium in choice}
(denoted SEC) if $x_{i}^{\ast }\in \bigcap\nolimits_{i\in N}C_{i}.$

\section{Equilibrium results}

This section provides a summary of different theorems concerning the
existence of equilibria in choice for games in choice form. In order to
underline the novelty and the importance of our work, we also must discuss
here the additional benefit of obtaining corollaries which state, under new
conditions, the existence of the weak Nash equilibria for games in normal
form, or the existence of the weak equilibria of qualitative games. To prove
our point, we will use continuous selection theorems or fixed point theorems
for the correspondences that we will form, considering upper sections of the
sets defining the game of the choice form. The advantage we have by doing
this, deserves a great prominence in the assessment of the new assumptions
which characterize the new statements. These results differ very much from
the ones obtained by Stefanescu, Ferrara and Stefanescu in \cite{SS} and
they link the recent research to the older approaches, regarding games in
normal form or qualitative games. In order to suggest priorities to the
reader, we keep the relationship between the main theorems and their
consequences on particular games. For a better understanding of the paper,
we recall all properties of the correspondences which will be used and the
tools of proofs. Our study gives a new perspective of unifying of different
approaches and results on the equilibrium concepts and the existence of
noncooperative theory of games. Finally, we note that we obtain the
existence of the strong equilibrium in choice for all situations considered
in this section, if we suppose, in addition, that $C_{i}(x_{-i})\neq
\emptyset $ for each $x_{-i}\in X_{-i}.\medskip $

Let $X$, $Y$ be topological spaces. We recall that the correspondence $%
T:X\rightarrow 2^{Y}$ has the \textit{local} \textit{intersection property }%
if $x\in X$ with $T(x)\neq \emptyset $ implies the existence of an open
neighborhood $V(x)$ of $x$ such that $\cap _{z\in V(x)}T(z)\neq \emptyset .$%
\medskip

To prove Theorem 1, we need the following lemma.

\begin{lemma}
(Wu, Shen, \cite{shen}). \textit{Let }$X$\textit{\ be a nonempty paracompact
subset of a Hausdorff topological space }$E$\textit{\ and }$Y$\textit{\ be a
nonempty subset of a Hausdorff topological vector space. Let }$S$\textit{,}$%
T:X\rightarrow 2^{Y}$\textit{\ be correspondence which verify:}
\end{lemma}

\textit{a) for each }$x\in X,$\textit{\ }$S(x)\neq \emptyset $\textit{\ and
co}$S(x)\subset T\left( x\right) ;$

\textit{b) }$S$\textit{\ has the local intersection property.}

\textit{Then, }$T$\textit{\ has a continuous selection.\medskip }

Our main result cites conditions which ensure the existence of equilibria in
choice for a game in choice form in the lack of convexity of the upper
sections of the sets $C_{i}$. The framework for our general next theorem
consists of Hausdorff topological vector spaces. The proof is based on an
argument that implies the above lemma.

\begin{theorem}
Let $((X_{i})_{i\in N},(C_{i})_{i\in N})$ be a \textit{game in choice form.
Assume that, for each }$i\in N,$ the following conditions are fulfilled:
\end{theorem}

\textit{a)} $X_{i}$ \textit{is a nonempty, convex and compact set in a} 
\textit{Hausdorff topological vector space }$E_{i}$;

\textit{b)} \textit{there exists a nonempty subset }$D_{i}$\textit{\ of }$%
C_{i}$ \textit{such that }$W_{i}=\{x_{-i}\in X_{-i}:D_{i}(x_{-i})\neq
\emptyset \}$\textit{\ is closed and }$D_{i}(x_{-i})\neq \emptyset $\textit{%
\ if only if} $C_{i}(x_{-i})\neq \emptyset ;$

\textit{c)} \textit{if} $D_{i}(x_{-i})\neq \emptyset ,$ \textit{there exists
an open neighborhood }$V(x_{-i})$\textit{\ of }$x_{-i}$\textit{\ such that }$%
\cap _{z_{-i}\in V(x_{-i})}D_{i}(z_{-i})\neq \emptyset ;$

\textit{d)} $D_{i}(x_{-i})$ \textit{is convex or empty for each} $x_{-i}\in
X_{-i}.$

\textit{Then, the game admits equilibria in choice.}

\begin{proof}
For each $i\in N,$ let us define the correspondences $S_{i},T_{i}:X_{-i}%
\rightarrow 2^{X_{i}},$ by

$S_{i}(x_{-i})=\left\{ 
\begin{array}{c}
\text{co}(\bigcup\nolimits_{\{y_{-i}:D_{i}(y_{-i})\neq \emptyset
\}}D_{i}(y_{-i}))\text{ if }x_{-i}\notin W_{i}; \\ 
D_{i}(x_{-i})\text{ if }x_{-i}\in W_{i}%
\end{array}%
\right. $ and

$T_{i}(x_{-i})=\left\{ 
\begin{array}{c}
\text{co}(\bigcup\nolimits_{\{y_{-i}:C_{i}(y_{-i})\neq \emptyset
\}}C_{i}(y_{-i}))\text{ if }x_{-i}\notin W_{i}; \\ 
C_{i}(x_{-i})\text{ if }x_{-i}\in W_{i}.%
\end{array}%
\right. $

We call $T_{i}$ the correspondence of the upper sections of the sets $C_{i}.$

The correspondence $S_{i}$ has nonempty and convex values and $%
S_{i}(x_{-i})\subset T_{i}(x_{-i})$ for each $x_{-i}\in X_{-i}.$

Assumption c) implies that $S_{i\mid (X_{-i}\cap W_{i})}$\textit{\ }has the
local intersection property. If $x_{-i}\in ^{C}W_{i},$ then, $S_{i}(x_{-i})=$%
co$(\bigcup\nolimits_{\{y_{-i}:D_{i}(y_{-i})\neq \emptyset
\}}D_{i}(x_{-i}))\neq \emptyset $. According to b), there exists an open
neighborhood\textit{\ }$V_{i}(x_{-i})$\textit{\ }of\textit{\ }$x_{-i}$%
\textit{\ }such that $D_{i}(z_{-i})=\emptyset $ for each $z_{-i}\in
V_{i}(x_{-i}).$ Then, $S_{i}(z_{-i})=$co$(\bigcup\nolimits_{%
\{y_{-i}:D_{i}(y_{-i})\neq \emptyset \}}D_{i}(y_{-i}))$ for each $z_{-i}\in
V_{i}(x_{-i})$ and, $\cap _{z_{-i}\in V_{i}(x_{-i})}S_{i}(z_{-i})=$co$%
(\bigcup\nolimits_{\{y_{-i}:D_{i}(y_{-i})\neq \emptyset
\}}D_{i}(y_{-i}))\neq \emptyset .$ It follows that $S_{i}:X_{-i}\rightarrow
2^{X_{i}}$ has the local intersection property.

The Wu-Shen Lemma implies that $T_{i}$\textit{\ }has a continuous selection $%
f_{i}:X_{-i}\rightarrow X_{i}.$

Let $f:X\rightarrow X$ be defined by $f(x)=\prod\nolimits_{i\in
N}f_{i}(x_{-i})$ for each $x\in X.$ The function $f$ is continuous, and,
according to the Brouwer fixed point Theorem, there exists $x^{\ast }\in X$
such that $f(x^{\ast })=x^{\ast }.$ Hence, $x^{\ast }\in
\prod\nolimits_{i\in N}T_{i}(x_{-i}^{\ast })$ and obviously, $x_{i}^{\ast
}\in T_{i}(x_{-i}^{\ast })$ for each $i\in N.$ Suppose that $(x_{-i}^{\ast
},X_{i})\cap C_{i}\neq \emptyset ,$ for some $i\in N.$ Then, $%
C_{i}(x_{-i}^{\ast })\neq \emptyset $ and $x_{i}^{\ast }\in
C_{i}(x_{-i}^{\ast }),$ which implies $x^{\ast }=(x_{-i}^{\ast },x_{i}^{\ast
})\in C_{i}.\medskip $
\end{proof}

As a corollary, we obtain sufficient conditions for a game in normal form to
admit weak Nash equilibria. The main assumption is new in literature and it
refers to the existence of a best reply for each player $i$ to the common
strategy of the other players, which lies in open intervals of the product
space $X_{-i}.$ The meaning is that each player $i$ can remain stable in the
choice of his own best strategy in the situation that the decisions of the
opponents can vary slightly in any manner profitable to themselves.

\begin{corollary}
Let $((X_{i})_{i\in N},(u_{i})_{i\in N})$ be a \textit{game in normal form.
Assume that, for each }$i\in N,$ the following conditions are fulfilled:
\end{corollary}

\textit{a)} $X_{i}$ \textit{is a nonempty, convex and compact set in a} 
\textit{Hausdorff topological vector space }$E_{i}$;

\textit{b)} \textit{the set }$\{x\in X:u_{i}(x)\geq u_{i}(x_{-i},y_{i})$ for
each $y_{i}\in X_{i}\}$ is nonempty;

\textit{c)} $W_{i}=\{x_{-i}\in X_{-i}:B_{i}(x_{-i})\neq \emptyset \}$\textit{%
\ is closed}$;$

\textit{d)} \textit{if} $B_{i}(x_{-i})\neq \emptyset ,$ \textit{there exists
an open neighborhood }$V(x_{-i})$\textit{\ of }$x_{-i}$\textit{\ so that }$%
\cap _{z_{-i}\in V(x_{-i})}B_{i}(z_{-i})\neq \emptyset ;$

\textit{e)} $B_{i}(x_{-i})$ \textit{is convex or empty for each} $x_{-i}\in
X_{-i}.$

\textit{Then, the game admits weak Nash equilibria.\medskip }

The following corollary concerns the existence of the weak equilibria\textit{%
\ }for qualitative games.

\begin{corollary}
Let $((X_{i})_{i\in N},(P_{i})_{i\in N})$ be a qualitative \textit{game.
Assume that, for each }$i\in N,$ the following conditions are fulfilled:
\end{corollary}

\textit{a)} $X_{i}$ \textit{is a nonempty, convex and compact set in a} 
\textit{Hausdorff topological vector space }$E_{i}$;

\textit{b)} \textit{the set }$\{x\in X:P_{i}(x)=\emptyset \}$ \textit{is
nonempty and }$W_{i}=\{x_{-i}\in X_{-i}:\exists x_{i}\in X_{i}$ \textit{such
that} $P_{i}(x_{-i},x_{i})=\emptyset \}$\textit{\ is closed}$;$

\textit{c)} for each $x_{-i}\in X_{-i},$ \textit{if} \textit{there exists }$%
x_{i}\in X_{i}$\textit{\ such that }$P_{i}(x)=\emptyset ,$\textit{\ then,
there exist an open neighborhood }$V(x_{-i})$\textit{\ of }$x_{-i}$\textit{\
and }$z_{i}\in X_{i}$ \textit{such that }$P_{i}(z_{-i},z_{i})=\emptyset $
for each $z_{-i}\in V_{i}(x_{-i});$

\textit{d)} $\{x_{i}\in X_{i}:P_{i}(x_{-i},x_{i})=\emptyset \}$ \textit{is
convex or empty for each} $x_{-i}\in X_{-i}.$

\textit{Then, the game admits weak equilibria.\medskip }

A new approach to the existence of equilibria in choice implies the concept
of weakly convex graph for a correspondence, proposed by X. Ding and Y He in 
\cite{126}. Firstly, we recall the definition.

\begin{definition}
\textit{(see }\cite{126}).\textit{\ }Let $X$ be a nonempty convex subset of
a topological vector space $E$ \ and $Y$ be a nonempty subset of $E$.\textit{%
\ }The correspondence\textit{\ }$T:X\longrightarrow 2^{Y}$ is said to have 
\emph{weakly convex graph} (in short, it is a WCG correspondence) if for
each $n\in \mathbb{N}$ and for each finite set $\{x_{1},x_{2},...,x_{n}\}%
\subset X$, there exists $y_{i}\in T(x_{i})$, $(i=1,2,...,n)$, such that
\end{definition}

\begin{center}
co$(\{(x_{1},y_{1}),(x_{2},y_{2}),...,(x_{n},y_{n})\})\subset $Gr$(T).$
\end{center}

We note that either the graph Gr$(T)$ is convex, or $\ \tbigcap \{T(x):x\in
X\}\neq \emptyset $, then $T$ has a weakly convex graph.

It is obvious that\textit{\ }a WCG\ correspondence may have nonempty values
and may not be convex-valued.

\begin{example}
$T:[0,2]\rightarrow \lbrack 0,2],$ $T(x)=\left\{ 
\begin{array}{c}
\lbrack 0,\frac{1}{2}]\cup \lbrack \frac{3}{2},2]\text{ \ \ \ \ \ if }x=1,
\\ 
\lbrack 0,2-\frac{1}{2}x]\text{ \ if }x\in \lbrack 0,2]\backprime \{1\}%
\end{array}%
\right. $ is a WCG correspondence (since $\cap \{T(x):x\in \lbrack
0,2]\}=\{0\}\neq \emptyset $), but $T(1)$ is not convex and Gr$(T)$ is not
convex either.
\end{example}

The theorem below is a continuous selection theorem for correspondences
having weakly convex graph.

\begin{lemma}
(Patriche, \cite{85}). \textit{Let }$Y$\textit{\ be a nonempty subset of a
topological vector space }$E$\textit{\ and }$K$\textit{\ be a }$(n-1)$%
\textit{- dimensional simplex in a topological vector space }$F$\textit{.
Let }$T:K\rightarrow 2^{Y}$\textit{\ be a WCG correspondence. Then, }$T$%
\textit{\ has a continuous selection on }$K$\textit{.\medskip }
\end{lemma}

The following lemma guarantees the existence of a fixed point for a product
of lower semi-continuous correspondences. It will be useful for proving our
second result on equilibria in choice.

\begin{lemma}
(Wu, \cite{120}) \textit{Let }$I$\textit{\ be an index set. For each }$i\in
I,$\textit{\ let }$X_{i}$\textit{\ be a nonempty convex subset of a
Hausdorff locally convex topological vector space }$E_{i}$\textit{, }$D_{i}$%
\textit{\ a nonempty compact metrizable subset of }$X_{i}$\textit{\ and }$%
S_{i},T_{i}:X\rightarrow 2^{D_{i}}$\textit{\ two correspondences with the
following conditions:}
\end{lemma}

(1) \textit{for each }$x\in X$, clco$S_{i}(x)\subset T_{i}(x)$ and $%
S_{i}(x)\neq \emptyset $\textit{. }

\textit{(2) }$S_{i}$\textit{\ is lower semi-continuous.}

\textit{Then, there exists} $x^{\ast }=\underset{i\in I}{\tprod }x_{i}^{\ast
}\in D=\underset{i\in I}{\tprod }D_{i}$ \textit{\ such that }$x_{i}^{\ast
}\in T_{i}(x^{\ast })$ for each $i\in I.\medskip $

The assumption that each correspondence of the upper sections of the sets $%
C_{i}$ has a selection which is an \textit{WCG }correspondence also assures
the existence of the equilibria in choice. The following theorem presents
precisely this.

\begin{theorem}
Let $((X_{i})_{i\in N},(C_{i})_{i\in N})$ be a \textit{game in choice form.
Assume that, for each }$i\in N,$ the following conditions are fulfilled:
\end{theorem}

\textit{a)} $X_{i}$\textit{\ is a nonempty convex compact metrizable set in
a locally convex space }$E_{i}$ \textit{and }$C_{i}$ \textit{is nonempty;}

\textit{b)} for each $x_{i}\in X_{i},$\textit{\ }$W_{i}=\{x_{-i}\in
X_{-i}:C_{i}(x_{-i})\neq \emptyset \}$ \textit{is} \textit{a }$(n_{i}-1)$%
\textit{-dimensional simplex in }$X;$

c) \textit{there exists a WCG correspondence }$S_{i}:W_{i}\rightarrow
2^{X_{i}}$\textit{\ such that }$S_{i}(x_{-i})\subset C_{i}(x_{-i})$ for each 
$x_{-i}\in W_{i};$

\textit{d)} $C_{i}(x_{-i})$ \textit{is convex or empty for each} $x_{-i}\in
X_{-i}.$

\textit{Then, the game admits equilibria in choice.}

\textit{Proof. }Let be $i\in I.$ From the assumption (c) and the selection
Lemma 3, it follows that there exists a continuous function $%
f_{i}:W_{i}\rightarrow X_{i}$ so that for each $x_{-i}\in W_{i}$, $%
f_{i}(x_{-i})\in S_{i}(x_{-i})\subset C_{i}(x_{-i}).$

Let us define the correspondence $T_{i}:X_{-i}\rightarrow 2^{X_{i}}$, by

$T_{i}(x_{-i}):=\left\{ 
\begin{array}{c}
\{f_{i}(x_{-i})\}\text{ \ \ \ \ \ \ \ \ \ \ \ \ \ if \ \ \ \ \ \ \ \ \ \ \ }%
x_{-i}\in W_{i}; \\ 
\text{co}(\bigcup\nolimits_{\{y_{-i}:C_{i}(y_{-i})\neq \emptyset
\}}C_{i}(y_{-i}))\text{ if }x_{-i}\notin W_{i}.%
\end{array}%
\right. $

We appeal to Lemma 1 to conclude that $T_{i}$ is lower semicontinuous on $X$%
. Tychonoff's Theorem implies that $X$ is compact.

According to Wu's fixed-point Theorem, applied for the correspondences $%
S_{i}=T_{i}$ and $T_{i}:X\rightarrow 2^{X_{i}},$ there exists $\mathit{x}%
^{\ast }\in X$ such that for each $i\in I$, $\mathit{x}_{i}^{\ast }\in T_{i}(%
\mathit{x}^{\ast })$. Suppose that $(x_{-i}^{\ast },X_{i})\cap C_{i}\neq
\emptyset ,$ for some $i\in N.$ Then, $C_{i}(x_{-i}^{\ast })\neq \emptyset $
and $x_{i}^{\ast }\in C_{i}(x_{-i}^{\ast }),$ which implies $x^{\ast
}=(x_{-i}^{\ast },x_{i}^{\ast })\in C_{i}.\medskip $

\begin{remark}
We can obtain two corollaries to the above theorem, if we replace assumption
c) with:
\end{remark}

c') \textit{there exists a correspondence }$S_{i}:W_{i}\rightarrow 2^{X_{i}}$%
\textit{\ such that}$\mathit{\ }S_{i}$ \textit{has a convex graph and }$%
S_{i}(x_{-i})\subset C_{i}(x_{-i})$ for each $x_{-i}\in W_{i};$

or

c") \textit{there exists a correspondence }$S_{i}:W_{i}\rightarrow 2^{X_{i}}$%
\textit{\ with closed values such that}$\mathit{\ }S_{i}$ \textit{has the
property that for any \ finite set }$\{x_{-i}^{1},x_{-i}^{2},...x_{-i}^{m}\}%
\subset X,$\textit{\ }$\tbigcap_{j=1}^{m}S_{i}(x_{-i}^{j})\neq \emptyset $%
\textit{\ and }$S_{i}(x_{-i})\subset C_{i}(x_{-i})$ for each $x_{-i}\in
W_{i};$

\textit{Proof}. In the first case, since a correspondence with a convex
graph is a WCG one, it follows that $S_{i}$ verifies Assumption c) from
Theorem 2, then we can apply this theorem.

In the second case, $X$ is a compact space and for each $i\in I$ the closed
sets $S_{i}(x_{-i}),$ $x_{-i}\in X_{-i}$ have the finite intersection
property, then $\tbigcap \{S_{i}(x_{-i}):x_{-i}\in X_{-i}\}\neq \emptyset $.
It follows that $S_{i}$ is a WCG correspondence and the conclusion comes
from Theorem 2.\medskip

As in the first case, we obtain the following corollaries concerning the
existence of the weak Nash equilibria for games in normal form and,
respectively, of the weak equilibria for qualitative games.

\begin{corollary}
\textit{Let }$((X_{i})_{i\in N},(u_{i})_{i\in N})$\textit{\ be a game in
normal form. Assume that, for each }$i\in N,$\textit{\ the following
conditions are fulfilled:}
\end{corollary}

\textit{a) }$X_{i}$\textit{\ is a nonempty, convex and compact set in a
Hausdorff topological vector space }$E_{i}$\textit{;}

\textit{b) the set }$\{x\in X:u_{i}(x)\geq u_{i}(x_{-i},y_{i})$\textit{\ for
each }$y_{i}\in X_{i}\}$\textit{\ is nonempty}$;$

\textit{c) for each }$x_{i}\in X_{i},$\textit{\ }$W_{i}=\{x_{-i}\in
X_{-i}:B_{i}(x_{-i})\neq \emptyset \}$\textit{\ is a }$(n_{i}-1)$\textit{%
-dimensional simplex in }$X;$

\textit{d) there exists a WCG correspondence }$S_{i}:W_{i}\rightarrow
2^{X_{i}}$\textit{\ such that }$S_{i}(x_{-i})\subset B_{i}(x_{-i})$\textit{\
for each }$x_{-i}\in W_{i};$

\textit{e) }$B_{i}(x_{-i})$\textit{\ is convex or empty for each }$x_{-i}\in
X_{-i}.$

\textit{Then, the game admits weak Nash equilibria.}

\begin{corollary}
\textit{Let }$((X_{i})_{i\in N},(P_{i})_{i\in N})$\textit{\ be a qualitative
game. Assume that, for each }$i\in N,$\textit{\ the following conditions are
fulfilled:}
\end{corollary}

\textit{a) }$X_{i}$\textit{\ is a nonempty, convex and compact set in a
Hausdorff topological vector space }$E_{i}$\textit{;}

\textit{b) the set }$\{x\in X:P_{i}(x)=\emptyset \}$\textit{\ is nonempty; }

\textit{c) for each }$x_{i}\in X_{i},$\textit{\ }$W_{i}=\{x_{-i}\in
X_{-i}:\exists x_{i}\in X_{i}$\textit{\ such that }$P_{i}(x_{-i},x_{i})=%
\emptyset \}$\textit{\ is a }$(n_{i}-1)$\textit{-dimensional simplex in }$X;$

\textit{d) there exists a WCG correspondence }$S_{i}:W_{i}\rightarrow
2^{X_{i}}$\textit{\ such that }$S_{i}(x_{-i})\subset \{x_{i}\in X_{i}$%
\textit{\ such that }$P_{i}(x_{-i},x_{i})=\emptyset \}$\textit{\ for each }$%
x_{-i}\in W_{i};$

\textit{e) }$\{x_{i}\in X_{i}:P_{i}(x_{-i},x_{i})=\emptyset \}$\textit{\ is
convex or empty for each }$x_{-i}\in X_{-i}.$

\textit{Then, the game admits weak equilibria.\medskip }

Now, we present a continuous selection lemma on Banach spaces which was
proposed by Yuan \cite{134}.

\begin{lemma}
(see \cite{134}). \textit{Let }$X$\textit{\ be a paracompact space, }$Y$%
\textit{\ be a Banach space and }$T:X\rightarrow 2^{Y}$\textit{\ be a lower
semicontinuous correspondences with closed convex values. Let }$%
S:X\rightarrow 2^{Y}$\textit{\ be a correspondence whose graph is open in }$%
X\times Y$\textit{\ such that }$T(x)\cap S(x)\neq \emptyset $\textit{\ for
each }$x\in X.$\textit{\ Then, there exists a continuous function }$%
f:X\rightarrow Y$\textit{\ such that }$f(x)\in $co$(T(x)\cap S(x))$\textit{\
for each }$x\in X$\textit{.\medskip }
\end{lemma}

The previous lemma leads us to the enunciation of Theorem 3, which gives new
conditions under which the equilibria in choice exist.

\begin{theorem}
\textit{Let }$((X_{i})_{i\in N},(C_{i})_{i\in N})$\textit{\ be a game in
choice form. Assume that, for each }$i\in N,$\textit{\ the following
conditions are fulfilled:}
\end{theorem}

\textit{a) }$X_{i}$\textit{\ is a nonempty, convex and compact set in a
Banach space }$E_{i};$\textit{\ }

\textit{b) }$C_{i}$\textit{\ is nonempty and open;}

\textit{c) the set }$W_{i}=\{x_{-i}\in X_{-i}:C_{i}(x_{-i})\neq \emptyset \}$%
\textit{\ is nonempty and open}$;$

\textit{d) there exists a lower semicontinuous correspondence }$%
S_{i}:W_{i}\rightarrow 2^{X_{i}}$\textit{\ with closed convex values such
that }$C_{i}(x_{-i})\cap S_{i}(x_{-i})\neq \emptyset $\textit{\ for each }$%
x_{-i}\in X_{-i}$\textit{;}

\textit{e) }$C_{i}(x_{-i})$\textit{\ is convex or empty for each }$x_{-i}\in
X_{-i}.$

\textit{Then, the game admits equilibria in choice.}

\begin{proof}
Let be $i\in I.$ From the assumption d) and the above lemma, it follows that
there exists a continuous function $f_{i}:W_{i}\rightarrow X_{i}$ such that $%
f_{i}(x_{-i})\in C_{i}(x_{-i})\cap S_{i}(x_{-i})$ for each $x_{-i}\in W_{i}$.

Let us define the correspondence $T_{i}:X_{-i}\rightarrow 2^{X_{i}}$, by

$T_{i}(x_{-i}):=\left\{ 
\begin{array}{c}
\{f_{i}(x_{-i})\}\text{ \ \ \ \ \ \ \ \ \ \ \ \ if \ \ \ \ \ \ \ \ \ \ \ \ }%
x_{-i}\in W_{i}; \\ 
\text{co}(\bigcup\nolimits_{\{y_{-i}:C_{i}(y_{-i})\neq \emptyset
\}}C_{i}(y_{-i}))\text{ if }x_{-i}\notin W_{i}.%
\end{array}%
\right. $

We appeal to Lemma 1 to conclude that $T_{i}$ is upper semicontinuous on $X$.

According to Kakutani's fixed-point Theorem, there exists $\mathit{x}^{\ast
}\in X$ such that for each $i\in I$, $\mathit{x}_{i}^{\ast }\in T_{i}(%
\mathit{x}_{-i}^{\ast })$. Suppose that $(x_{-i}^{\ast },X_{i})\cap
C_{i}\neq \emptyset ,$ for some $i\in N.$ Then, $C_{i}(x_{-i}^{\ast })\neq
\emptyset $ and $x_{i}^{\ast }=f_{i}(x_{-i}^{\ast })\in C_{i}(x_{-i}^{\ast
}),$ which implies $x^{\ast }=(x_{-i}^{\ast },x_{i}^{\ast })\in
C_{i}.\medskip $
\end{proof}

As corollaries of Theorem 3, we obtain the following results which assume
the lower semicontinuity of the involved correspondences. We \ note that
these results are different from the old ones obtained in literature so far.

\begin{corollary}
\textit{Let }$((X_{i})_{i\in N},(u_{i})_{i\in N})$\textit{\ be a game in
normal form. Assume that, for each }$i\in N,$\textit{\ the following
conditions are fulfilled:}
\end{corollary}

\textit{a) }$X_{i}$\textit{\ is a nonempty, convex and compact set in a
Banach space }$E_{i}$\textit{;}

\textit{b) the set }$\{x\in X:u_{i}(x)\geq u_{i}(x_{-i},y_{i})$\textit{\ for
each }$y_{i}\in X_{i}\}$\textit{\ is nonempty and open}$;$

\textit{c) the set }$W_{i}=\{x_{-i}\in X_{-i}:B_{i}(x_{-i})\neq \emptyset \}$%
\textit{\ is nonempty and open}$;$

\textit{d) there exists a lower semicontinuous correspondence }$%
S_{i}:X_{-i}\rightarrow 2^{X_{i}}$\textit{\ with closed convex values such
that }$S_{i}(x_{-i})\cap B_{i}(x_{-i})\neq \emptyset $\textit{\ for each }$%
x_{-i}\in X_{-i}$\textit{;}

\textit{e) }$B_{i}(x_{-i})$\textit{\ is convex or empty for each }$x_{-i}\in
X_{-i}.$

\textit{Then, the game admits weak Nash equilibria.}

\begin{corollary}
\textit{Let }$((X_{i})_{i\in N},(P_{i})_{i\in N})$\textit{\ be a qualitative
game. Assume that, for each }$i\in N,$\textit{\ the following conditions are
fulfilled:}
\end{corollary}

\textit{a) }$X_{i}$\textit{\ is a nonempty, convex and compact set in a
Banach space }$E_{i}$\textit{;}

\textit{b) the set }$\{x\in X:P_{i}(x)=\emptyset \}$\textit{\ is nonempty
and open; }

\textit{c) the set }$W_{i}=\{x_{-i}\in X_{-i}:$ $\mathit{\exists }$ $%
x_{i}\in X_{i}$ \textit{such that} $P_{i}(x_{-i},x_{i})=\emptyset \}$\textit{%
\ is nonempty and open}$;$

\textit{d) there exists a lower semicontinuous correspondence }$%
S_{i}:W_{i}\rightarrow 2^{X_{i}}$\textit{\ with closed convex values such
that }$S_{i}(x_{-i})\cap \{x_{i}\in X_{i}:P_{i}(x_{-i},x_{i})=\emptyset
\}\neq \emptyset $\textit{\ for each }$x_{-i}\in X_{-i}$\textit{;}

\textit{e) }$\{x_{i}\in X_{i}:P_{i}(x_{-i},x_{i})=\emptyset \}$\textit{\ is
convex or empty for each }$x_{-i}\in X_{-i}.$

\textit{Then, the game admits weak equilibria.\medskip }

Further, we will prove the existence of the equilibrium in choice under new
conditions. The proof we will provide explicitly relies on \ lemmas
concerning the fixed points for the correspondences we will construct based
on the upper sections of the sets $(C_{i})_{i\in N}.$

First, we recall the following definition.

If $X$ is a nonempty set and $Y$ is a topological space, the correspondence $%
T:X\rightarrow 2^{Y}$ is said to be \textit{transfer open-valued} \cite{tian}
if for any $(x,y)\in X\times Y$ with $y\in T(x),$ there exists an $x^{\prime
}\in X$ such that $y\in $int$T(x^{\prime }).\medskip $

Further, we present the following useful statement about the transfer
open-valued correspondences (Proposition 1 in \cite{lin2001}).

\begin{lemma}
Let $Y$ be a nonempty set, $X$ be a topological space and $S:X\rightarrow
2^{Y}$ be a correspondence. The following assertions are equivalent:
\end{lemma}

\textit{a) }$S^{-1}:Y\rightarrow 2^{X}$\textit{\ is transfer open-valued and 
}$S$\textit{\ has nonempty values;}

\textit{b) }$X=\bigcup\nolimits_{y\in Y}$int$S^{-1}(y).\medskip $

In this context, Ansari and Yao proved in \cite{an} a fixed point
result.\medskip

\begin{lemma}
(Ansari, Yao \cite{an}). \textit{Let }$X$\textit{\ be a compact convex
subset of a Hausdorff topological vector space. Let }$S:X\rightarrow 2^{X}$%
\textit{\ be a correspondence with nonempty convex values.\ If }$X=\cup \{$%
int$_{X}S^{-1}(y):y\in X\}$ (or, $S^{-1}:Y\rightarrow 2^{X}$\textit{\ is
transfer open-valued)}$,$ \textit{then, }$S$\textit{\ has fixed
points.\medskip }
\end{lemma}

We will apply the previous lemma in order to prove the existence of the
equilibria in choice for games in choice form.

\begin{theorem}
Let $((X_{i})_{i\in N},(C_{i})_{i\in N})$ be a \textit{game in choice form.
Assume that, for each }$i\in N,$ the following conditions are fulfilled:
\end{theorem}

\textit{a)} $X_{i}$ \textit{is a nonempty, convex and compact set in a} 
\textit{Hausdorff topological vector space }$E_{i}$ \textit{and }$C_{i}$ 
\textit{is nonempty;}

\textit{b)} $X=\bigcup\nolimits_{y\in X}\{$int$_{X}\bigcap\nolimits_{i\in
N}(^{C}W_{i}\cup \{x_{-i}\in X_{-i}:x_{-i}\in C_{i}(y_{i})\})\}$\textit{,
where }$W_{i}=\{x_{-i}\in X_{-i}:C_{i}(x_{-i})\neq \emptyset \};$

\textit{c)} $C_{i}(x_{-i})$ \textit{is convex or empty for each} $x_{-i}\in
X_{-i}.$

\textit{Then, the game admits equilibria in choice.}

\begin{proof}
For each $i\in N,$ let us define the correspondence $T_{i}:X_{-i}\rightarrow
2^{X_{i}},$ by

$T_{i}(x_{-i})=\left\{ 
\begin{array}{c}
\text{co}(\bigcup\nolimits_{\{y_{-i}:C_{i}(y_{-i})\neq \emptyset
\}}C_{i}(y_{-i}))\text{ if }x_{i};\notin W_{i} \\ 
C_{i}(x_{-i})\text{ \ \ \ \ \ \ \ \ \ \ \ \ \ \ if \ \ \ \ \ \ \ \ \ \ \ \ \
\ }x_{-i}\in W_{i}.%
\end{array}%
\right. $

The correspondence $S_{i}$ has nonempty and convex values$.$

Let $T:X\rightarrow 2^{X}$ be defined by $T(x)=\prod\nolimits_{i\in
N}T_{i}(x_{-i})$ for each $x\in X.$

The correspondence $S$ also has nonempty and convex values$.$

If $y\in X,$ then $T^{-1}(y)=\bigcap\nolimits_{i\in N}\{x\in X:y_{i}\in
T_{i}(x_{-i})\}=\bigcap\nolimits_{i\in N}(^{C}W_{i}\cup \{x_{-i}\in
X_{-i}:x_{-i}\in C_{i}(y_{i})\}).$

$X=\bigcup\nolimits_{y\in X}$int$_{X}T^{-1}(y)$, according to assumption b).

We can apply the Ansari and Yao Lemma and we obtain that there exists $%
x^{\ast }\in X$ such that $x^{\ast }\in T(x^{\ast }).$ Obviously, $%
x_{i}^{\ast }\in T_{i}(x_{-i}^{\ast })$ for each $i\in N.$ Suppose that $%
(x_{-i}^{\ast },X_{i})\cap C_{i}\neq \emptyset ,$ for some $i\in N.$ Then, $%
C_{i}(x_{-i}^{\ast })\neq \emptyset $ and $x_{i}^{\ast }\in
C_{i}(x_{-i}^{\ast }),$ which implies $x^{\ast }=(x_{-i}^{\ast },x_{i}^{\ast
})\in C_{i}.$
\end{proof}

\begin{remark}
According to Lemma 7, we can replace condition b) in Theorem 4 with
\end{remark}

\textit{b')} if $x_{i}\in C_{i}(y_{-i}),$\textit{\ then, there exists }$%
z_{i}\in X_{i}$ \textit{such that }$y_{-i}\in $\textit{int}$%
_{X_{-i}}\{x_{-i}\in X_{-i}:z_{i}\in C_{i}(x_{-i})\}$\textit{\ and the set }$%
W_{i}=\{x_{-i}\in X_{-i}:C_{i}(x_{-i})\neq \emptyset \}$\textit{\ is
closed.\medskip }

A new result involving the equilibria in choice\textit{\ }will naturally
follow directly from Theorem 4.

\begin{theorem}
Let $((X_{i})_{i\in N},(C_{i})_{i\in N})$ be a \textit{game in choice form.
Assume that, for each }$i\in N,$ the following conditions are fulfilled:
\end{theorem}

\textit{a)} $X_{i}$ \textit{is a nonempty, convex and compact set in a} 
\textit{Hausdorff topological vector space }$E_{i}$ \textit{and }$C_{i}$ 
\textit{is nonempty;}

\textit{b)} \textit{for each }$x_{i}\in X_{i},$\textit{\ }$\{x_{-i}\in
X_{-i}:x_{i}\in C_{i}(x_{-i})\}\cup ^{C}W_{i}$\textit{\ is open, where }$%
W_{i}=\{x_{-i}\in X_{-i}:C_{i}(x_{-i})\neq \emptyset \};$

\textit{c)} $C_{i}(x_{-i})$ \textit{is convex or empty for each} $x_{-i}\in
X_{-i}.$

\textit{Then, the game admits equilibria in choice.\medskip }

Another proof of the above theorem appeals to Yannelis and Prabhakar'
continuous selection lemma \cite{127} \ applied for the correspondences $%
T_{i}:X_{-i}\rightarrow 2^{X_{i}}$ defined by

$T_{i}(x_{-i})=\left\{ 
\begin{array}{c}
\text{co}(\bigcup\nolimits_{\{y_{-i}:C_{i}(y_{-i})\neq \emptyset
\}}C_{i}(y_{-i}))\text{ if }x_{-i}\notin W_{i}; \\ 
C_{i}(x_{-i})\text{ \ \ \ \ \ \ \ \ \ \ \ \ \ \ if \ \ \ \ \ \ \ \ \ \ \ \ \
\ }x_{-i}\in W_{i}.%
\end{array}%
\right. $ $\ $for each $i\in N.$

We present below the lemma.\medskip

\begin{lemma}
(Yannelis and Prabhakar, \cite{127}). \textit{Let }$X$\textit{\ be a
paracompact Hausdorff topological space and }$Y$\textit{\ be a Hausdorff
topological vector space. Let }$T:X\rightarrow 2^{Y}$\textit{\ be a
correspondence with nonempty convex values\ and for each }$y\in Y$\textit{, }%
$T^{-1}(y)$\textit{\ is open in }$X$\textit{. Then, }$T$\textit{\ has a
continuous selection that is, there exists a continuous function }$%
f:X\rightarrow Y$\textit{\ so that }$f(x)\in T(x)$\textit{\ for each }$x\in
X $\textit{.\medskip }
\end{lemma}

For the proof, we note that for each $i\in N,$ the correspondence $T_{i}$
has nonempty and convex values and if $x_{i}\in X_{i},$ then $%
T_{i}^{-1}(x_{i})=^{C}W_{i}\cup \{x_{-i}\in X_{-i}:x_{i}\in C_{i}(x_{-i})\}$
is an open set, according to assumption b).

We can apply the Yannelis and Prabhakar Lemma and we obtain that there
exists $f_{i}:X_{-i}\rightarrow X_{i}$, a continuous selection of $T_{i}.$
Let $f:X\rightarrow X$ be defined by $f(x)=\prod\nolimits_{i\in
N}f_{i}(x_{-i})$ for each $x\in X.$ The function $f$ is continuous, and,
according to the Brouwer fixed point Theorem, there exists $x^{\ast }\in X$
such that $f(x^{\ast })=x^{\ast }.$ Hence, $x^{\ast }\in
\prod\nolimits_{i\in N}T_{i}(x_{-i}^{\ast })$ and obviously, $x_{i}^{\ast
}\in T_{i}(x_{-i}^{\ast })$ for each $i\in N.$ Suppose that $(x_{-i}^{\ast
},X_{i})\cap C_{i}\neq \emptyset ,$ for some $i\in N.$ Then, $%
C_{i}(x_{-i}^{\ast })\neq \emptyset $ and $x_{i}^{\ast }\in
C_{i}(x_{-i}^{\ast }),$ which implies $x^{\ast }=(x_{-i}^{\ast },x_{i}^{\ast
})\in C_{i}.$

\begin{remark}
If for each $x_{i}\in X_{i},$\ $\{x_{-i}\in X_{-i}:x_{i}\in C_{i}(x_{-i})\}$%
\ is open, and the set $W_{i}=\{x_{-i}\in X_{-i}:C_{i}(x_{-i})\neq \emptyset
\}$ is closed, then the above Theorem holds.\medskip
\end{remark}

Corollary 7 is mainly obtained by verifying an assumption concerning the
union of all lower sections of the best reply correspondences for a game in
normal form\textit{.}

\begin{corollary}
\textit{Let }$((X_{i})_{i\in N},(u_{i})_{i\in N})$\textit{\ be a game in
normal form. Assume that, for each }$i\in N,$\textit{\ the following
conditions are fulfilled:}
\end{corollary}

\textit{a) }$X_{i}$\textit{\ is a nonempty, convex and compact set in a
Hausdorff topological vector space }$E_{i}$\textit{;}

\textit{b) the set }$\{x\in X:u_{i}(x)\geq u_{i}(x_{-i},y_{i})$\textit{\ for
each }$y_{i}\in X_{i}\}$\textit{\ is nonempty}$;$

\textit{c) }$X=\bigcup\nolimits_{y\in X}\{$int$_{X}\bigcap\nolimits_{i\in
N}(^{C}W_{i}\cup B_{i}^{-1}(y_{i}))\}$\textit{, where }$W_{i}=\{x_{-i}\in
X_{-i}:B_{i}(x_{-i})\neq \emptyset \};$

\textit{d) }$B_{i}(x_{-i})$\textit{\ is convex or empty for each }$x_{-i}\in
X_{-i}.$

\textit{Then, the game admits weak Nash equilibria.}

\begin{remark}
In Corollary 7, condition c) can also be replaces with
\end{remark}

c') \textit{the best reply correspondence }$B_{i}:X_{-i}\rightarrow
2^{X_{i}} $\textit{\ is transfer open valued and the set} $W_{i}=\{x_{-i}\in
X_{-i}:B_{i}(x_{-i})\neq \emptyset \}$\textit{\ is closed}$.$

\medskip A new statement can be deduced explicitly from Corollary 7.

\begin{corollary}
\textit{Let }$((X_{i})_{i\in N},(u_{i})_{i\in N})$\textit{\ be a game in
normal form. Assume that, for each }$i\in N,$\textit{\ the following
conditions are fulfilled:}
\end{corollary}

\textit{a) }$X_{i}$\textit{\ is a nonempty, convex and compact set in a
Hausdorff topological vector space }$E_{i}$\textit{;}

\textit{b) the set }$\{x\in X:u_{i}(x)\geq u_{i}(x_{-i},y_{i})$\textit{\ for
each }$y_{i}\in X_{i}\}$\textit{\ is nonempty}$;$

\textit{c) }$B_{i}^{-1}(x_{i})\cup \{x_{-i}\in
X_{-i}:B_{i}(x_{-i})=\emptyset \}$\textit{\ is open for each }$x_{i}\in
X_{i};$

\textit{d) }$B_{i}(x_{-i})$\textit{\ is convex or empty for each }$x_{-i}\in
X_{-i}.$

\textit{Then, the game admits weak Nash equilibria.\medskip }

The following results refer to the existence of weak equilibria for the
qualitative games. They are consequences of Theorem 5.

\begin{corollary}
\textit{Let }$((X_{i})_{i\in N},(P_{i})_{i\in N})$\textit{\ be a qualitative
game. Assume that, for each }$i\in N,$\textit{\ the following conditions are
fulfilled:}
\end{corollary}

\textit{a) }$X_{i}$\textit{\ is a nonempty, convex and compact set in a
Hausdorff topological vector space }$E_{i}$\textit{;}

\textit{b) the set }$\{x\in X:P_{i}(x)=\emptyset \}$\textit{\ is nonempty; }

\textit{c) }$X=\bigcup\nolimits_{y\in X}\{$int$_{X}\bigcap\nolimits_{i\in
N}(^{C}W_{i}\cup \{x_{-i}\in X_{-i}:$\textit{\ }$P_{i}(x_{-i},y_{i})=%
\emptyset \})\},$\textit{\ where }$W_{i}=\{x_{-i}\in X_{-i}:\exists x_{i}\in
X_{i}$ \textit{such that }$P_{i}(x_{-i},x_{i})=\emptyset \};$\textit{\ }

\textit{d) }$\{x_{i}\in X_{i}:P_{i}(x_{-i},x_{i})=\emptyset \}$\textit{\ is
convex or empty for each }$x_{-i}\in X_{-i}.$

\textit{Then, the game admits weak equilibria.}

\begin{remark}
In Corollary 9, condition c) can also be replaced with
\end{remark}

c")\textit{if} $P_{i}(y_{-i},x_{i})=\emptyset ,$\textit{\ then, there exists 
}$z_{i}\in X_{i}$ \textit{such that }$y_{-i}\in $\textit{int}$%
_{Y_{-i}}\{x_{-i}\in X_{-i}:P_{i}(x_{-i},z_{i})=\emptyset \}$\textit{\ and
the set }$W_{i}=\{x_{-i}\in X_{-i}:\exists x_{i}\in X_{i}$\textit{\ such
that }$P_{i}(x_{-i},x_{i})=\emptyset \}$\textit{\ is open.}

\begin{corollary}
\textit{Let }$((X_{i})_{i\in N},(P_{i})_{i\in N})$\textit{\ be a qualitative
game. Assume that, for each }$i\in N,$\textit{\ the following conditions are
fulfilled:}
\end{corollary}

\textit{a) }$X_{i}$\textit{\ is a nonempty, convex and compact set in a
Hausdorff topological vector space }$E_{i}$\textit{;}

\textit{b) the set }$\{x\in X:P_{i}(x)=\emptyset \}$\textit{\ is nonempty; }

\textit{c) }$\{x_{-i}\in X_{-i}:\nexists x_{i}\in X_{i}$\textit{\ such that }%
$P_{i}(x_{-i},x_{i})=\emptyset \}\cup \{x_{-i}\in X_{-i}:$\textit{\ }$%
P_{i}(x_{-i},x_{i})=\emptyset \}$\textit{\ is open for each }$x_{i}\in
X_{i}; $\textit{\ }

\textit{d) }$\{x_{i}\in X_{i}:P_{i}(x_{-i},x_{i})=\emptyset \}$\textit{\ is
convex or empty for each }$x_{-i}\in X_{-i}.$

\textit{Then, the game admits weak equilibria.\medskip }

A very great importance in the fixed point theory has Tarafdar's fixed point
Theorem, which we present below.

\begin{lemma}
(Tarafdar, \cite{114}). \textit{Let }$\{X_{i}\}_{i\in I}$\textit{\ be a
family of nonempty compact convex sets, each in a topological vector space }$%
E_{i},$\textit{\ where }$\mathit{I}$\textit{\ is an index set. Let }$%
X=\prod\limits_{i\in I}X_{i}.$\textit{\ For each }$i\in I,$\textit{\ let }$%
S_{i}:X\rightarrow 2^{X_{i}}$\textit{\ be a correspondence such that}
\end{lemma}

a) \textit{for each }$x\in X,$\textit{\ }$S_{i}(x)$\textit{\ is a nonempty,
convex subset of }$X_{i};$

b) \textit{for each }$x_{i}\in X_{i},$\textit{\ }$S_{i}^{-1}(x_{i})$\textit{%
\ contains a relatively open subset }$O_{x_{i}}$\textit{\ of }$X$\textit{\
such that}

$\cup _{x_{i}\in X_{i}}O_{x_{i}}=X$\textit{\ }$(O_{x_{i}}$\textit{\ may be
empty for some }$x_{i}).$

\textit{Then, there exists a point }$x\in X$\textit{\ \ such that }$x\in
S(x)=\prod\limits_{i\in I}S_{i}(x)$\textit{, that is, }$x_{i}\in S_{i}(x)$%
\textit{\ for each }$i\in I,$\textit{\ where }$x_{i}$\textit{\ is the
projection of }$x$\textit{\ onto }$X_{i}$\textit{\ for each }$i\in I.$

\medskip By using Lemma 9, we establish Theorem 6, which is a slightly
different version of Theorem 4.

\begin{theorem}
Let $((X_{i})_{i\in N},(C_{i})_{i\in N})$ be a \textit{game in choice form.
Assume that, for each }$i\in N,$ the following conditions are fulfilled:
\end{theorem}

\textit{a)} $X_{i}$ \textit{is a nonempty, convex and compact set in a} 
\textit{Hausdorff topological vector space }$E_{i}$ \textit{and }$C_{i}$ 
\textit{is nonempty;}

\textit{b)} for each $x_{i}\in X_{i},$\textit{\ }$\{x_{-i}\in
X_{-i}:x_{i}\in C_{i}(x_{-i})\}$ \textit{contains a relatively open subset }$%
O_{x_{i}}$\textit{\ of }$X_{-i}$,\textit{\ such that}

$\cup _{x_{i}\in X_{i}}O_{x_{i}}=X_{-i}$\textit{\ }$(O_{x_{i}}$\textit{\ may
be empty for some }$x_{i});$

\textit{c)} $C_{i}(x_{-i})$ \textit{is convex or empty for each} $x_{-i}\in
X_{-i}.$

\textit{Then, the game admits equilibria in choice.}

\begin{proof}
For each $i\in N,$ let us define the correspondence $T_{i}:X_{-i}\rightarrow
2^{X_{i}},$ by

$T_{i}(x_{-i})=\left\{ 
\begin{array}{c}
\text{co}(\bigcup\nolimits_{\{y_{-i}:C_{i}(y_{-i})\neq \emptyset
\}}C_{i}(y_{-i}))\text{ if }x_{-i}\notin W_{i}; \\ 
C_{i}(x_{-i})\text{ \ \ \ \ \ \ \ \ \ \ \ \ \ \ if \ \ \ \ \ \ \ \ \ \ \ \ \
\ }x_{-i}\in W_{i},%
\end{array}%
\right. ,$ \textit{where }$W_{i}=\{x_{-i}\in X_{-i}:C_{i}(x_{-i})\neq
\emptyset \}.$

The correspondence $T_{i}$ has nonempty and convex values$.$

If $x_{i}\in X_{i},$ then $T_{i}^{-1}(x_{i})=^{C}W_{i}\cup \{x_{-i}\in
X_{-i}:x_{i}\in C_{i}(x_{-i})\}.$

According to assumption b), for each $x_{i}\in X_{i},$\ $T_{i}^{-1}(x_{i})$\
contains a relatively open subset $O_{x_{i}}$\ of $X$\ such that

$\cup _{x_{i}\in X_{i}}O_{x_{i}}=X_{-i}$\ $(O_{x_{i}}$\ may be empty for
some $x_{i}).$ Then, $X=\cup _{x_{i}\in X_{i}}(O_{x_{i}}\times X_{i}).$ For
each $i\in N,$ let us define the correspondence $S_{i}:X\rightarrow
2^{X_{i}},$ by $S_{i}(x)=T_{i}(x_{-i}).$ If $x_{i}\in X_{i},$ then $%
S_{i}^{-1}(x_{i})=T_{i}^{-1}(x_{i})\times X_{i}.$ We denote $%
U_{x_{i}}=O_{x_{i}}\times X_{i}$ and we obtain that for each $x_{i}\in
X_{i}, $\textit{\ }$S_{i}^{-1}(x_{i})$ contains a relatively open subset $%
U_{x_{i}}$\ of $X$\ such that $\cup _{x_{i}\in X_{i}}U_{x_{i}}=X$\ $%
(U_{x_{i}}$\ may be empty for some $x_{i}).$

We can apply the previous lemma and we obtain that there exists $x^{\ast
}\in X$ such that $x_{i}^{\ast }\in S_{i}(x^{\ast })=T_{i}(x_{-i}^{\ast })$
for each $i\in N.$ Suppose that $(x_{-i}^{\ast },X_{i})\cap C_{i}\neq
\emptyset ,$ for some $i\in N.$ Then, $C_{i}(x_{-i}^{\ast })\neq \emptyset $
and $x_{i}^{\ast }\in C_{i}(x_{-i}^{\ast }),$ which implies $x^{\ast
}=(x_{-i}^{\ast },x_{i}^{\ast })\in C_{i}.\medskip $
\end{proof}

Now, we get the following corollaries from the previous result.

\begin{corollary}
\textit{Let }$((X_{i})_{i\in N},(u_{i})_{i\in N})$\textit{\ be a game in
normal form. Assume that, for each }$i\in N,$\textit{\ the following
conditions are fulfilled:}
\end{corollary}

\textit{a) }$X_{i}$\textit{\ is a nonempty, convex and compact set in a
Hausdorff topological vector space }$E_{i}$\textit{;}

\textit{b) the set }$\{x\in X:u_{i}(x)\geq u_{i}(x_{-i},y_{i})$\textit{\ for
each }$y_{i}\in X_{i}\}$\textit{\ is nonempty}$;$

\textit{c)} \textit{for each} $x_{i}\in X_{i},$\textit{\ }$B_{i}^{-1}(x_{i})$
\textit{contains a relatively open subset }$O_{x_{i}}$\textit{\ of }$X$%
\textit{\ such that}

$\cup _{x_{i}\in X_{i}}O_{x_{i}}=X_{-i}$\textit{\ }$(O_{x_{i}}$\textit{\ may
be empty for some }$x_{i});$

\textit{d) }$B_{i}(x_{-i})=\{x_{i}\in X_{i}:u_{i}(x)\geq u_{i}(x_{-i},y_{i})$
for each $y_{i}\in X_{i}\}$\textit{\ is convex or empty for each }$x_{-i}\in
X_{-i}.$

\textit{Then, the game admits weak Nash equilibria.}

\begin{corollary}
\textit{Let }$((X_{i})_{i\in N},(P_{i})_{i\in N})$\textit{\ be a qualitative
game. Assume that, for each }$i\in N,$\textit{\ the following conditions are
fulfilled:}
\end{corollary}

\textit{a) }$X_{i}$\textit{\ is a nonempty, convex and compact set in a
Hausdorff topological vector space }$E_{i}$\textit{;}

\textit{b) the set }$\{x\in X:P_{i}(x)=\emptyset \}$\textit{\ is nonempty; }

\textit{c)} \textit{for each} $x_{i}\in X_{i},$\textit{\ }$\{x_{-i}\in
X_{-i}:P_{i}(x_{-i},x_{i})=\emptyset \}$ \textit{contains a relatively open
subset }$O_{x_{i}}$\textit{\ of }$X$\textit{\ such that}

$\cup _{x_{i}\in X_{i}}O_{x_{i}}=X_{-i}$\textit{\ }$(O_{x_{i}}$\textit{\ may
be empty for some }$x_{i});$

\textit{d) }$\{x_{i}\in X_{i}:P_{i}(x_{-i},x_{i})=\emptyset \}$\textit{\ is
convex or empty for each }$x_{-i}\in X_{-i}.$

\textit{Then, the game admits weak equilibria.}

\begin{remark}
In a particular case, we can weaken condition b) of Theorem 6 by condition
b'):
\end{remark}

b') for each $x_{i}\in X_{i},$\textit{\ }$\{x_{-i}\in X_{-i}:x_{i}\in
C_{i}(x_{-i})\}=O_{x_{i}}$ is an open subset of $X$\ such that

$\cup _{x_{i}\in X_{i}}O_{x_{i}}=X_{-i}$\ $(O_{x_{i}}$\ may be empty for
some $x_{i}).$ According to Lemma 7, this condition is equivalent with the
fact that the correspondence $T_{i}^{-1}:X_{i}\rightarrow 2^{X_{-i}}$\textit{%
\ }is transfer open-valued and $T_{i}$\ has nonempty values, where\textit{\ }%
$T_{i}:X_{-i}\rightarrow 2^{X_{i}}$ is defined by $%
T_{i}(x_{-i})=C_{i}(x_{-i})$ for each $x_{-i}\in X_{-i}$\textit{. }

In this case, we obtain the following theorem concerning the existence of
the strong equilibrium in choice.

\begin{theorem}
Let $((X_{i})_{i\in N},(C_{i})_{i\in N})$ be a \textit{game in choice form.
Assume that, for each }$i\in N,$ the following conditions are fulfilled:
\end{theorem}

\textit{a)} $X_{i}$ \textit{is a nonempty, convex and compact set in a} 
\textit{Hausdorff topological vector space }$E_{i}$ \textit{and }$C_{i}$ 
\textit{is nonempty;}

\textit{b)} for each $x_{i}\in X_{i},$\textit{\ }$\{x_{-i}\in
X_{-i}:x_{i}\in C_{i}(x_{-i})\}$ $=O_{x_{i}}$ \textit{is an open subset of }$%
X$\textit{\ such that}

$\cup _{x_{i}\in X_{i}}O_{x_{i}}=X_{-i}$\textit{\ }$(O_{x_{i}}$\textit{\ may
be empty for some }$x_{i});$

\textit{c)} $C_{i}(x_{-i})$ \textit{is convex or empty for each} $x_{-i}\in
X_{-i}.$

\textit{Then, the game admits strong equilibria in choice.\medskip }

\section{Concluding remarks}

Our study is a new perspective unifying different approaches and results on
the equilibrium concepts and the existence of noncooperative theory of
games. We have proposed to the reader a synthesis of theorems and
consequences which state, under new conditions, the existence of the
equilibrium for games in choice form, in normal form and also for
qualitative games. Our approach differs essentially from the one of
Stefanescu, Ferrara and Stefanescu, \ who proposed the new concept of game
in choice form and the corresponding equilibrium in choice (2012). A further
research may consist of the integration of all research instruments and
perspectives. The advantage of using these new ideas is that they are more
systematic and can cover more general situations. This paper reflects the
integrity of this kind of thinking and can reopen the problem of the
equilibrium existence under new conditions.

\end{document}